\documentclass{amsproc}
\usepackage{amssymb}
\newtheorem{theorem}{Theorem}[section]

\newtheorem{prp}[theorem]{Proposition}
\newtheorem{crl}[theorem]{Corollary}
\newtheorem{cnj}[theorem]{Conjecture}

\theoremstyle{definition}
\newtheorem{dfn}[theorem]{Definition}
\newtheorem{example}[theorem]{Example}

\theoremstyle{remark}
\newtheorem{remark}[theorem]{Remark}

\numberwithin{equation}{section}

\newcommand{\bin}[2]{\left({#1}\atop{#2}\right)}

\newcommand{\prj}[1]{{\rm Proj}\left(#1^{[*]}\right)}
\begin{document}

\title{Combinatorics and Geometry\\of Higher Level Weyl Modules}

\author{B. Feigin}
\address{Landau institute for Theoretical Physics, Chernogolovka
142432, Russia}
\email{feigin@mccme.ru}
\thanks{The first author was supported in part by Grants RFBR-05-01-01007, 
SS~2044-2003-2 and INTAS~03-51-3350}

\author{A.N.~Kirillov}
\address{RIMS, Kyoto University, Kyoto 606-8502, Japan}
\email{kirillov@kurims.kyoto-u.ac.jp}
\thanks{The second author was supported in part by Grand-in-Aid-15340006 for
Scientific Research  from Japan Society for Promotion of Science (JSPS)}

\author{S. Loktev}
\address{Department of Mathematics, Kyoto University, Kyoto 606-8502, Japan}
\email{loktev@math.kyoto-u.ac.jp}
\thanks{The third  author was supported in part by Grants CRDF 
RM1-2545-MO-03, PICS~2094, RFBR-04-01-00702, Russian Federation President Grant MK-3419.2005.1 and Grant-in-Aid for
Exploratory Research no.~16654003 from the Ministry of Education,
Culture, Sports, Science and Technology, Japan.}

\date{\today}



\begin{abstract}
A  higher level analog of Weyl modules over multi-variable currents is proposed. 
It is shown that the sum of their dual spaces
form a commutative algebra. The structure of these modules and the geometry
 of the projective spectrum of this algebra is
studied for the currents of dimension one and two. Along the way we prove
some particular cases of the conjectures in \cite{FL1} and propose
a generalization of the notion of parking function representations.
\end{abstract}

\maketitle

\def \a  {{\mathfrak a}}
\def \CC {{\mathbb C}}
\def \g {{\mathfrak g}}
\def \b {{\mathfrak b}}
\def \h {{\mathfrak h}}
\def \n {{\mathfrak n}}
\def \oo {{\mathcal O}}
\def \U {{\mathbf U}}
\def \p {{\mathfrak p}}
\def \z {{\mathcal Z}}
\def \A  {{\mathcal A}}

\def \gr {{\rm gr}}
\def \Gr {{\rm Gr}}
\def \gl {{\mathfrak gl}}
\def \sh {{\mathcal P}}
\def \Sh {{\mathrm Sh}}
\def \fm {{\mathcal H}}
\def \fu {{\Gamma}}
\def \hilb {{\mathcal H}}
\def \grc {\gr_C}
\def \nv {{\mathcal N}_w}
\def \ms {{\mathcal M}}
\def \ch {{\rm ch}}

\def \mat {{\rm Mat}}

\section*{Introduction}

Let us start from some geometrical background.  
The Borel--Weyl theorem provides a construction of finite--dimensional representations
of a simple Lie algebra in spaces of sections of line bundles on the corresponding flag
variety. This theorem can be generalized for affine algebras, but the affine flag varieties
are infinite--dimensional and it is not so convenient to work with them. One approach is to
treat this variety as a limit of finite--dimensional Schubert subvarieties. And the
space of  sections of those line bundles on these subvarieties are known as {\em Demazure modules}.

The affine Demazure modules are well studied (see \cite{KMOTU}, \cite{FLitt}, \cite{S}), but 
there was a lack of construction for such 
modules, which are not bases on infinite--dimensional representation theory. Now it appears
(\cite{CL}, \cite{FLitt2}) that for the simply-laced case the level one Demazure modules are isomorphic to the classical
Weyl modules introduced in \cite{CP}. 
Then the higher level Demazure modules can be constructed 
in a standard way (See Section 1.1 and Section 2).

Another construction, based on the fusion product introduced in \cite{FL1}, is
related to the tensor product structure as a module over constant currents (see \cite{CP}
for Weyl modules and \cite{FLitt} for Demazure modules). Here we show that in the simply-laced case 
the fusion product of irreducible representations
 produces the Demazure modules (see \cite{CL}, \cite{FLitt2} for the level one case).

In \cite{FL2} an analog of Weyl modules for multi--dimensional currents was introduced. It can be
considered as a candidate for level one Demazure modules over the multi--dimensional current 
algebra. Here we construct higher level Weyl modules that pretend to be Demazure modules
of arbitrary level. In particular, their
dual spaces form a commutative algebra, whose spectrum  can be considered as a multi--dimensional
analog of a corresponding Schubert variety.

In the case of double affine (toroidal) algebras things become exciting. Recall that the toroidal
algebra is the universal central extension of $\g \otimes \CC[x,x^{-1}, y, y^{-1}]$.
An analog of level one integrable representations (and their quantum version)
was introduced in \cite{VV}, \cite{STU}. 

Weyl modules are also studied for that case in \cite{FL2}, \cite{FL3}. For $\g = sl_r$
its structure was established for the weights, proportional to the weight of vector
representation, and a conjecture for other weights was proposed. 
Since the construction in \cite{FL3} is pretty similar to given in \cite{VV}, \cite{STU},
we are convinced that 
the Weyl modules over two-dimensional polynomial ring are isomorphic to $\g \otimes \CC[x, y]$--submodules of 
that toroidal modules. So the Weyl modules can be considered as an analog
of level one Demazure modules also in toroidal settings.
Moreover, in \cite{FL3} an action of the universal central extension of
$\g \otimes \CC[x, x^{-1},y]$ on the limits of Weyl modules is proposed. We expect that these limits
are just the restriction of the corresponding toroidal modules to $\g \otimes \CC[x, x^{-1},y]$.

Recall that in \cite{FL3} we construct $gl_r \otimes \CC[x, y]$--modules from cyclic modules over
$gl_r \otimes B_N$, where $B_N$ is the associative algebra of upper--triangular $N\times N$ matrices.
Note that $gl_r \otimes B_N$ is isomorphic to a parabolic subalgebra $\p \subset gl_{Nr}$, 
so we can obtain such representations from the spaces of sections of $\p$--equivariant bundles.
Namely, suppose we have a line bundle on a closure of a  $\p$--orbit. Then the space, dual
to its section, is a cyclic representation of  $\p$, so we can produce a
representation of $gl_r \otimes \CC[x, y]$ from it.
And our conjecture is that the higher level Weyl modules can be obtained 
in this way from Schubert varieties.

We start Section~1 from some general notion and constructions. Namely, 
we introduce higher level cyclic modules and a structure of commutative algebra on
their dual spaces. Then we discuss geometry of the spectrum of this algebra and propose
some useful examples. At the end of the section we recall the definition of Weyl modules.

In Section~2 we discuss in more detail 
higher level Weyl modules and fusion modules over one--dimensional currents
with values in a simply-laced simple Lie algebra. Here we relate
them each to other as well as to Demazure modules. In $sl_r$ case we construct each
higher level Weyl module as a fusion module as well as a Demazure module.

In Section~3 we proceed to two--dimensional currents. Here we relate the higher 
level Weyl modules
to the coordinate ring of the usual Schubert varieties in a Grassmann variety, using
the deformation of Weyl modules proposed in \cite{FL3}. Also
we introduce a higher level generalization of the parking function notion.

{\bf Acknowledgments} We are grateful to V.~Chari, G. Fourier, P. Littelmann and 
M.~Okado for useful discussions. The third
author is on leave from Independent University of Moscow and Institute for Theoretical
and Experimental Physics.

\section{Generalities}

Recall that a module, generated by a single vector, is called a {\em cyclic} module
and this vector is called a {\em cyclic} vector.

\subsection{Higher level modules}

\begin{dfn}
Let $W$ be a cyclic module over a Lie algebra $\a$
generated by a fixed cyclic vector $w$. Introduce the module $W^{[k]}$
as the submodule of $W^{\otimes k}$ generated by $w^{\otimes k}$
\end{dfn}

\begin{prp}\label{prop_qu}
Suppose $W_0$ is a quotient of $W$, fix the image $w_0$ of $w$
as the cyclic vector in $W_0$. Then $W_0^{[k]}$ is a quotient of $W^{[k]}$.
\end{prp}

\begin{proof}
We have the natural map of modules $W^{\otimes k} \to W_0^{\otimes k}$ sending
$w^{\otimes k}$ to $w_0^{\otimes k}$. Then the image of $W^{[k]}$ is the submodule of
$W_0^{\otimes k}$, generated by $w_0^{\otimes k}$, that is $W_0^{[k]}$.
\end{proof}

Now suppose that the action of $\a$ on $W$ is extended to an action of
a corresponding connected Lie group $\A$.

\begin{prp}\label{prop_lg}
We have $W^{[k]}$ is the $\A$--submodule in $W^{\otimes k}$ generated by $w^{\otimes k}$.
\end{prp}

\begin{proof}
Proceeding to the tangent spaces, we obtain that this $\A$--submodule contains  $W^{[k]}$.
 On the other hand, since $\A$
is generated by the image of the exponent map, we have the inverse inclusion.
\end{proof}

\subsection{Multiplication}

Note that we have the natural inclusion 
$$m_{k_1,k_2} : W^{[k_1+k_2]} \hookrightarrow W^{[k_1]} \otimes W^{[k_2]}$$ 
of subspaces in $W^{\otimes (k_1+k_2)}$.

\begin{prp}
Let $W^{[*]} = \bigoplus_{k\ge 0} \left( W^{[k]}\right)^*$, where $W^{[0]} = \CC$. The map 
$$m^* = \bigoplus_{k_1,k_2} m_{k_1,k_2}^* : 
W^{[*]} \otimes W^{[*]} \to W^{[*]}$$
defines a structure of commutative algebra on $W^{[*]}$.
\end{prp}

\begin{proof}
Let $m_{k_1,k_2,k_3}: W^{[k_1+k_2+k_3]} \hookrightarrow
 W^{[k_1]} \otimes W^{[k_2]} \otimes W^{[k_3]}$
be the similar inclusion.
Then we have
$$m_{k_1,k_2+k_3}^* (1\otimes m_{k_2,k_3}^*) = m_{k_1,k_2,k_3}^* =
 m_{k_1+k_2,k_3}^* (m_{k_1,k_2}^* \otimes 1).$$
So $m^* (1 \otimes m^*) = m^* (m^* \otimes 1)$, that is associativity.
\end{proof}

We can also view the algebra $W^{[*]}$ as the coordinate ring
of the projective variety $\prj{W}$.

Suppose that $\a$ is the Lie algebra of a connected Lie group $\A$ acting on $W$.
 Let $\oo(k)$ be the line bundles on $P(W)$
formed by homogeneous polynomials of degree $k$.  
Consider the $\A$--orbit of $\CC w$ in the projective space $P(W)$. 
Note that the closure $\nv$ of this orbit is algebraic.

\begin{prp}
We have $\prj{W} \cong  \nv$ and $W^{[k]} \cong \Gamma(\nv, \oo(k))^*$ as $\A$--modules.
\end{prp}

\begin{proof}
First we have the natural isomorphism $ev: S^k (W) \to \Gamma(P(W), \oo(k))^*$. Here 
$ev(u^{\otimes k})$ 
evaluate the sections at the point $\CC u \in P(W)$, and $S^k(W)$ is spanned by such vectors.

The space $\Gamma(\nv, \oo(k))$ is the quotient of $\Gamma(P(W), \oo(k))$ by 
the subspace of sections vanishing
at $\nv$. On the other hand, taking into account Proposition~\ref{prop_lg},
we have $W^{[k]}$ is the subspace of $S^k (W)$ spanned 
by $u^{\otimes k}$ with $\CC u \in \nv$ 
(we can proceed to the closure because any linear subspace is closed).
So they are dual each to other.

Concerning the multiplication, due to the action of $\a$ the maps $m_{k_1,k_2}$ are
 uniquely defined 
by the image of $w^{\otimes k}$. And for $\gamma_1 \in \oo(k_1)$, $\gamma_2 \in \oo(k_2)$ we have
$ev(w^{\otimes k}) (\gamma_1 \cdot \gamma_2) = ev(w^{\otimes k_1}) (\gamma_1)\, ev(w^{\otimes k_2}) (\gamma_2)$.
Since the multiplication of sections is also equivariant,
we have an isomorphism of algebras.
\end{proof}

\subsection{Examples}

\begin{example}
Let $\a = \g$ be a simple Lie algebra, let $V(\lambda)$ be the irreducible finite--dimensional representation with highest weight
$\lambda$. Fix the highest weight vector in $V(\lambda)$ as the cyclic one. 
Then $V(\lambda)^{[k]} \cong V(k\lambda)$ and 
$\prj{V(\lambda)}$ is known as a
{\em generalized flag variety} of the corresponding group $G$ 
(that is the usual flag variety if $\lambda$ is big enough).
\end{example}

\begin{example}\label{ex_sch}
Let $\a = \b \subset \g$ be the Borel subalgebra (that stabilizes the highest weight vector). For an element ${\bf w}$ of the
Weyl group  introduce the {\em extremal} vector
$v_{\bf w} \in V(\lambda)$ as the vector with the weight obtained from the highest weight by the action of ${\bf w}$. 
This vector is defined uniquely up to multiplication by a scalar.

Take $v_{\bf w}$ as the cyclic vector of $V_{v_{\bf w}}(\lambda) = U(\b) v_{\bf w}
 \subset V(\lambda)$.
Then we have $V_{v_{\bf w}}(\lambda)^{[k]} = V_{v_{\bf w}}(k\lambda)$ 
and $\prj{V_{v_{\bf w}}(\lambda)}$ is known as a {\em Schubert subvariety}
of the generalized flag variety.

Note that these varieties are well--defined and the same situation holds
for affine Lie algebras. 
\end{example}

\begin{example}\label{ex_gr}
Let $\g = \gl_r$, Take $\lambda = \omega_n$, that is the $n$--th fundamental weight, 
so $V(\lambda) = \wedge^n V$, where $V$ is the $r$--dimensional vector representation of
$\gl_r$. Then $\prj{V(\lambda)}$ is the Grassmann variety $
\Gr(n,r)$ formed by $n$--dimensional planes in $\CC^r$.
The restriction of $\oo(1)$ to $\Gr(n,r)$ 
is dual to the determinant bundle with stalks $\wedge^n P$ over each plane $P$.

Choose a basis $v_1, \dots, v_r$ in $V$, then fix the Borel subalgebra $\b$
of upper--triangular matrices mapping each $v_i$ to a linear combination of $v_j$ with $j\le i$.
Then extremal vectors in $\wedge^n V$ are just monomials of the form 
$v_\eta = v_{\eta_1} \wedge \dots \wedge v_{\eta_n}$, 
$\eta_1 < \dots < \eta_n$. The corresponding Schubert subvariety 
$\Sh_\eta = \prj{V_{v_\eta} (\omega_n)}$ 
consists of planes, whose intersection with
$\left< v_1, v_2, \dots v_{\eta_i}\right>$ has dimension at least $i$ for all $i=1\dots n$.
\end{example}

\subsection{Cyclic adjoint module}

Suppose we have an increasing
 filtration on the Lie algebra $\a$: $F^0 \a \subset F^1 \a  \subset \dots$. 
Then it can be extended to the filtration $F^0 U(\a) \subset F^1 U(\a) \subset \dots$.

For a cyclic module $W$ introduce the filtration and the adjoint graded space 
$$F^i_C W = (F^i U(\a))  w, \qquad \grc W = \bigoplus_i F^i_C W/  F^{i-1}_C W.
$$
Then $\grc W$ is a module over $\gr A$.

For $u \in W$ by $\overline{u}$ denote the corresponding vector in $\grc W$. Fix
$\overline{w}$ as a cyclic vector in $\grc W$.

\begin{prp}\label{prp_grd}
We have $(\grc W)^{[k]}$ is a quotient of $\grc W^{[k]}$.
\end{prp}  

\begin{proof}
First let us construct a map
$$\rho: \grc W^{[k]} \to (\grc W)^{\otimes k}.$$
We have 
$$F^i_C W^{[k]} \subset \sum_{i_1+ \dots+i_n = i} 
F^{i_1}_C  W \otimes \dots \otimes F^{i_n}_C  W.$$
 As
$$\left(F^{i_1}_C  W \otimes \dots \otimes F^{i_n}_C  W\right)\cap
\left(F^{i_1'}_C W \otimes \dots \otimes F^{i_n'}_C  W\right)
= F^{\min(i_1,i_1')}_C  W \otimes \dots \otimes F^{\min(i_n,i_n')}_C  W,$$
it gives a map 
$$F^i_C W^{[k]} \to \bigoplus_{i_1+ \dots+i_n = i} 
\grc^{i_1}  W \otimes \dots \otimes \grc^{i_n}  W,$$
where $\grc^{i}  W = F_C^i W /  F_C^{i-1} W$.
As the image of $F^{i-1}_C W^{[k]}$ under this map is zero, we obtain
the map $\rho$.

Note that $\rho(\overline{w^{\otimes k}}) = \overline{w}^{\otimes k}$, so the
image of $\rho$ is $(\grc W)^{[k]}$
and we have the statement of the proposition.
\end{proof}

\subsection{Weyl modules}

Let $\g$ be a reductive Lie algebra. Choose a Cartan and a Borel subalgebra 
$\h \subset \b \subset \g$.
In this paper we consider mainly $\a = \g \otimes \CC[x^1, \dots, x^d]$ and the following
class of modules.

\begin{dfn}
For a weight $\lambda: \b \to \h \to \CC$ 
let  $W^d(\lambda)$ be the maximal finite-dimensional module over 
$\g \otimes \CC[x^1, \dots, x^d]$ generated by $w_\lambda$ such that
\begin{equation}\label{hwc}
(\g \otimes P) w_\lambda = \lambda(g) P(0) w_\lambda \quad \mbox{for} \ g \in \b.
\end{equation}
\end{dfn}

By {\em maximal} we mean that any finite-dimensional module generated by $w_\lambda$ is
a quotient of $W^d(\lambda)$.

In \cite{FL2} it is shown that $W^d(\lambda)$ exists and that it is non-trivial for
a dominant $\lambda$. Also it is shown there that $W^d(\lambda)$ is graded as 
$\g \otimes \CC[x^1, \dots, x^d]$--module, that is we have
$$W^d(\lambda) = \bigoplus_{i_1, \dots, i_d \ge 0} W^d(\lambda)^{i_1, \dots, i_d},$$
where $W^d(\lambda)^{i_1, \dots, i_d}$ are $\g \otimes 1$--modules and
$\g \otimes (x^1)^{j_1} \dots (x^d)^{j_d}$ acts from
$W^d(\lambda)^{i_1, \dots, i_d}$ to $W^d(\lambda)^{i_1+j_1, \dots, i_d+j_d}$.

This grading can be extended in the usual way to the grading on the tensor product and therefore 
on the higher level Weyl modules.

Finally, let us introduce the notation for graded character of any graded module $W$ by
$$\ch_d W = \sum _{i_1, \dots, i_d \ge 0} t_1^{i_1} \dots t_d^{i_d} \cdot
\ch W^{i_1, \dots, i_d},$$
where $\ch$ denotes the usual character of $\g$--modules.

\section{One-dimensional case}

\def \wg {\widehat{\g}}

In this section let $\g$ be a simple simply-laced Lie algebra (that is of type A, D or E), 
let $\h \subset \g$ be a Cartan subalgebra.
By $R$ denote the set of roots and
by $Q$ denote the root lattice of $\g$.
Let $\omega_1, \dots \omega_r$ be the fundamental weights.

\subsection{Demazure modules}

Let $\wg = \g \otimes \CC [x^{-1},x] \oplus \CC c$ be the central extension of $\g \otimes \CC [x^{-1},x]$.
 Note that the restriction of the
central extension to $\g \otimes \CC[x]$ is trivial, so we have $\g \otimes \CC[x] \subset \wg$.
Let $L_k$ be the integrable level $k$ vacuum representation of $\wg$, let 
$w_k \in L_k$ be the highest weight vector. Recall that $c$ acts on $L_k$
by the scalar $k$ and $w_k$ is annihilated by
the subalgebra $\g \otimes \CC[x]$ (see \cite{K} for details about $\wg$ and $L_k$).

\def \i {\iota}

For a weight $\lambda$ let $\i_\lambda: \wg \to \wg$ be the automorphism mapping 
$g\otimes x^n$ to $g \otimes x^{n+\lambda(\alpha)}$ when $g$ belongs to
the root space 
$\g_\alpha$, $\alpha \in R$, so it stabilizes $h \otimes x^n$ when $h \in \h$, $n \ne 0$
and maps $h \otimes 1$ to $h \otimes 1 + \lambda(h) c$. 
Note that it coincides with the action of the corresponding translation from 
the extended affine Weyl group.
Let $\g[x] = \g \otimes \CC[x]$ and let $\g[x]_\lambda = \i_{\lambda}(\g \otimes \CC[x])$.
In other words, we have
$$\g[x]_\lambda = \h_\lambda  \oplus \h \otimes x\CC[x]\oplus
\bigoplus_{\alpha \in R}  \g_\alpha 
\otimes x^{\lambda(\alpha)} \CC[x] \subset \wg,$$
where $\h_\lambda$ consists of $h  \otimes 1 + \lambda(h) c$ for $h \in \h$.

\begin{dfn}
Introduce the {\em Demazure module} over  $\g[x]$ by
$$D(k,\lambda) = \i_{\lambda}^* \left(U(\g[x]_\lambda)w_k \right).$$
Namely, it is isomorphic to $U(\g[x]_\lambda)w_k$ as the vector space, where $\g[x]$ acts via
identification $\i_{\lambda}$ with $\g[x]_\lambda$.
\end{dfn}

Note that the classical highest weight of $D(k,\lambda)$ is $\lambda$, that is
$$(h\otimes 1) w_k = \lambda(h) w_k \quad \mbox{for} \ h \in \h.$$
When $\g$ 
is not simply-laced, this definition needs a modification and 
only a part of weights can be obtained as classical highest weights.

\begin{prp}
We have $D(1,\lambda)^{[k]} \cong D(k,\lambda)$.
\end{prp}

\begin{proof}
It follows from  $L_{1}^{[k]} \cong L_{k}$. And this is because $L_{1}^{[k]}$ is integrable 
and generated by the highest weight vector with the corresponding highest weight. 
\end{proof}

\begin{prp} (see also \cite{FF3}, \cite{S})
Suppose $\lambda \in Q$. Then $\prj{D(1,\lambda)}$ is a Schubert cell 
in the affine Grassmann variety for $\wg$.
\end{prp} 

\begin{proof}
For $\lambda \in Q$ the automorphism $\i_\lambda$ is the action of an element of the affine Weyl group, 
so we are in the situation of Example~\ref{ex_sch}.
\end{proof}

Concerning the dimension, one can use the following result.

\begin{theorem} \cite{FLitt}\label{thm_tns} For a dominant $\lambda$ we have
$$\dim D(k,\lambda) = \prod_{i=1}^r \left(\dim D(k,\omega_i)\right)^{\lambda_i}.$$ \qed
\end{theorem}

\begin{prp} For a dominant $\lambda$ we have
$D(k,\lambda)$ is a quotient of $W^1(\lambda)^{[k]}$.
\end{prp}

\begin{proof}
The module  $D(1,\lambda)$ is finite--dimensional and, since $\lambda(\alpha)\ge 0$
for a positive $\alpha$, the image of $w_1$ in $D(1,\lambda)$
satisfies~\eqref{hwc}, so $D(1,\lambda)$ is a quotient of  $W^1(\lambda)$.

As $D(k,\lambda)\cong D(1,\lambda)^{[k]}$, the proposition follows from 
Proposition~\ref{prop_qu}.
\end{proof}

\begin{cnj}\label{conj_d} {\rm (recently proved in \cite{FLitt2})}
For a dominant $\lambda$ we have
$$\dim D(1, \lambda) = \dim W^1(\lambda).$$
\end{cnj}

This conjecture is already proved in $sl_r$ case (\cite{CL}). Also note that we expect 
this equality only for a simply-laced $\g$.
Summarizing the statements above, we
have the following description of the higher level Weyl modules. 

\begin{prp}
In the case when Conjecture~\ref{conj_d} holds, for a dominant $\lambda$ we have
$W^1(\lambda)^{[k]}\cong  D(k,\lambda)$, therefore for  $\lambda \in Q$ 
we have $\prj{W^1(\lambda)}$  is a Schubert subvariety 
in the affine Grassmann variety.
\end{prp}

\subsection{Fusion}

\def \ac {{\g\otimes \CC[x]}}

For $z \in \CC^d$ let $\varphi(z)$ be
the automorphism of $\ac$ sending $x$ to $x+z$.
For a module $W$ over $\ac$ define the shifted module
$W(z) = \varphi(z)^* W$, so  $W(z) \cong W$ as a vector space and
the action is combined with $\varphi(z)$.

Let $W_1, \dots, W_n$ be cyclic modules over  $\ac$
 such that for a certain $N$ the subalgebra $\g \otimes x^N \CC[x]$ 
acts on each $W_i$ by zero. By $w_1, \dots , w_n$ denote their cyclic vectors.

\begin{prp} \cite{FL1}
If $z_i$ are pairwise distinct then the vector $w_1 \otimes \dots \otimes w_n$ is cyclic in 
$W_1(z_1) \otimes \dots \otimes W_n(z_n)$. 
\end{prp}

\begin{proof}
Let 
$$\g[N(z_1, \dots, z_n)] = \g \otimes \left(\CC[x]/ (x-z_1)^N \dots (x-z_n)^N \CC[x]
\right).$$

As the ideal $\g \otimes (x-z_1)^N \dots (x-z_n)^N \CC[x]$ acts on the tensor produce by zero,
$W_1(z_1) \otimes \dots \otimes W_n(z_n)$ is a module over $\g[N(z_1, \dots, z_n)]$.

Next note that we have the natural projections $p_i: \g[N(z_1, \dots, z_n)] \to\g[N(z_i)]$.  
Then it is known that the direct sum 
$$\oplus p_i: \g[N(z_1, \dots, z_n)] \to \bigoplus_{i=1}^n \g[N(z_i)]$$
is an isomorphism. 
Note that the preimage of each $\g[N(z_i)]$ belongs to the ideal
$\g \otimes (x-z_1)^N \dots (x-z_{i-1})^N(x-z_{i+1})^N\dots (x-z_n)^N \CC[x]$.
So we have
\begin{equation}\label{localize}
U( \g[N(z_1, \dots, z_n)]) (w_1 \otimes \dots \otimes w_n) = 
\left(U(\g[N(z_1)]\right) w_1 \otimes \dots \otimes \left(U(\g[N(z_n)]\right) w_n
\end{equation}
that is equal to  $W_1(z_1) \otimes \dots \otimes W_n(z_n)$.
\end{proof}

Introduce the filtration on $\g \otimes \CC[x]$ such that $F^i (\g \otimes \CC[x])$ consists
of $\g$-valued polynomials whose degree does not exceed $i$. Then $\gr (\g \otimes \CC[x]) =
\g \otimes \CC[x]$, so we can produce a graded module from any cyclic module.

\begin{dfn}\cite{FL1}
For given $z_1, \dots, z_n$ introduce the  $\g \otimes \CC[x]$-module
$$W_1 * \dots * W_n (z_1, \dots, z_n) = 
\grc
(W_1(z_1) \otimes \dots \otimes W_n(z_n)).$$
We call it the {\em fusion module}.
\end{dfn}

\begin{prp}\label{prop_fq}
For  any $\g \otimes \CC[x]$--modules $W_1, \dots, W_n$ we have the module
$\left(W_1 * \dots * W_n (z_1, \dots, z_n)\right)^{[k]}$
is a quotient of $W_1^{[k]} * \dots * W_n^{[k]} (z_1, \dots, z_n)$.
\end{prp}

\begin{proof}
Note that due to the isomorphism~\eqref{localize} we have
$$\left(W_1(z_1) \otimes \dots \otimes W_n(z_n)\right)^{[k]} \cong 
W_1^{[k]}(z_1) \otimes \dots \otimes W_n^{[k]}(z_n).$$
Then the statement follows from Proposition~\ref{prp_grd}.
\end{proof}

For any $\g$--module $V$ introduce the evaluation $\g \otimes \CC[x]$--module $V[0]$, 
isomorphic to $V$ as a vector space, where
$\g \otimes x\CC[x]$ acts by zero and $\g \otimes 1$ acts as $\g$ on $V$. For 
a set of $\g$--modules $V_1, \dots, V_n$ let us set
$$V_1 * \dots * V_n (z_1, \dots, z_n) = V_1[0] * \dots * V_n[0] (z_1, \dots, z_n).$$

\begin{cnj}\label{conj_v}
For any simple $\g$ and any weights $\lambda^1, \dots, \lambda^n$ we have
$$\left(V(\lambda^1) * \dots * V(\lambda^n) (z_1, \dots, z_n)\right)^{[k]} 
\cong V(k\lambda^1) * \dots * V(k\lambda^n) (z_1, \dots, z_n).$$
\end{cnj}

Note that $V(k\lambda) \cong V(\lambda)^{[k]}$, so by Proposition~\ref{prop_fq}
the left hand side in this Conjecture is a quotient of the right hand side. 

\subsection{Fusion of Weyl modules}

\begin{prp}\label{prp_qf} For any   $\lambda^1, \dots, \lambda^n$
we have $W^1(\lambda^1) * \dots * W^1(\lambda^n)(z_1, \dots, z_n)$ 
is a quotient of $W^1(\lambda^1+ \dots + \lambda^n)$.
\end{prp}

\begin{proof}
One can show that the cyclic vector of the left hand side satisfies~\eqref{hwc}. As this module is
finite--dimensional, the proposition follows from the definition of Weyl modules.
\end{proof}

Together with Proposition~\ref{prp_grd} it motivates the following conjecture.

\begin{cnj}\label{conj_f}
We have  
$$W^1(\lambda^1)^{[k]} * \dots * W^1(\lambda^n)^{[k]}(z_1, \dots, z_n)
 \cong W^1(\lambda^1+ \dots + \lambda^n)^{[k]},$$
in particular, the left hand side is independent on $z_1, \dots, z_n$.
\end{cnj}

\begin{theorem}
Conjecture~\ref{conj_d} implies Conjecture~\ref{conj_f}.
\end{theorem}

\begin{proof}

In the case $k=1$ Conjecture~\ref{conj_d} together with Theorem~\ref{thm_tns} implies
the equality of dimensions, so this case of Conjecture~\ref{conj_f} follows from
Proposition~\ref{prp_qf}.

Note that by Proposition~\ref{prop_fq}
the module
$W^1(\lambda^1)^{[k]}* \dots *   W^1(\lambda^n)^{[k]}(z_1, \dots, z_n)$
has the quotient
$$W^1(\lambda^1)* \dots *  W^1(\lambda^n)(z_1, \dots, z_n)^{[k]} \cong 
D(1,\lambda^1 + \dots +\lambda^n)^{[k]},$$
which is isomorphic to 
$D(k,\lambda^1 + \dots +\lambda^n)$ as well as to 
$W^1(\lambda^1 + \dots +\lambda^n)^{[k]}$.

By Theorem~\ref{thm_tns} we have
$$\dim D(k,\lambda^1 + \dots +\lambda^n) = \dim D(k, \lambda^1) \cdots \dim
D(k, \lambda^n),$$
that by assumption is equal to
$\dim W^1(\lambda^1)^{[k]}* \dots *   W^1(\lambda^n)^{[k]}(z_1, \dots, z_n)$.
So this quotient is the whole space and
we have the isomorphism proposed in Conjecture~\ref{conj_f}.
\end{proof}

\subsection{$gl_r$ case}

Now suppose $\g = \gl_r$.

\begin{theorem}\cite{CL} For $\g = \gl_r$ we have $W^1(\lambda) \cong D(1,\lambda)$.
\end{theorem}
So Conjecture~\ref{conj_d} and therefore Conjecture~\ref{conj_f} is already proved in
this case.

Note that for  $\g = \gl_r$ we have $W(\omega_i) \cong V(\omega_i)[0]$, where the right hand side
is the evaluation representation defined above.

\begin{crl}
We have 
$$ 
V(k\omega_{1})^{*\lambda_1}  * \dots * V(k\omega_{r-1})^{*\lambda_{r-1}}(z_1, \dots, z_{|\lambda|})
\cong D(k,\lambda),$$
in particular, the left hand side does not depend on $z_1, \dots, z_{|\lambda|}$.\qed
\end{crl}
\noindent So we proved a substantial case of 
Conjecture~1.8 from \cite{FL1}.
And we can deduce the similar case of Conjecture~\ref{conj_v}.

\begin{crl}
For $1\le i_1, \dots, i_n \le r-1$ we have
$$\left(V(l\omega_{i_1}) * \dots * V(l\omega_{i_n}) (z_1, \dots, z_n)\right)^{[k]}
\cong V(kl\omega_{i_1}) * \dots * V(kl\omega_{i_n}) (z_1, \dots, z_n).$$ \qed
\end{crl}

\begin{remark}
It seems that these two conjectures can be proved also for the set of weights
$l_1\omega_{i_1}, \dots, l_n\omega_{i_n}$ using the methods of \cite{FF2}, where
fusion product embedded into a direct sum of integrable modules. Then
the corresponding varieties $\prj{V(l_1\omega_{i_1}) * \dots * V(l_n\omega_{i_n}) 
(z_1, \dots, z_n)}$ 
coincide with the {\em generalized Schubert varieties}
introduced and described for $\g = sl_2$ in \cite{FF3}, \cite{F}.
\end{remark}

\begin{remark}
Note that due to the result of \cite{S} on Demazure modules
we also have a formula for the graded character 
$\ch_1 V(k\omega_{1})^{*\lambda_1}  * \dots * V(k\omega_{r-1})^{*\lambda_{r-1}}$ in terms of
parabolic Kostka polynomials as expected in \cite{FL1}.
\end{remark}

\section{Two-dimensional case}

Now let us consider the case $d=2$ and $\g = \g l_r$. Here we use the partition notation
for weights of $\g l_r$.

\subsection{Deformation of Weyl modules}
Let us recall a construction from \cite{FL3}.

The Lie algebra $\gl_r \otimes \CC[x^1,x^2]$ can be deformed into the Lie algebra
$\gl_r \otimes \CC\left< X,Y\right>$, where $\CC\left< X,Y\right>$ is the associative
algebra, generated by $X$ and $Y$ under the relation $YX-XY = X$. 
The algebra $\CC\left< X,Y\right>$ has the natural representation in
$\CC[t,t^{-1}]$, where $X$ acts by multiplication on  $t$ and $Y$ acts as $t \partial/
\partial t$.

By $V$ denote the $r$--dimensional vector representation of $\g l_r$. 
Let $v_1,\dots,v_r$ be the standard basis vectors in $V$.
For a partition $\xi$ introduce the $\gl_r \otimes \CC\left< X,Y\right>$--module
$V_\xi$ as the submodule of 
$$\bigwedge\nolimits^{|\xi|}\left(V \otimes  \CC[t,t^{-1}]/\CC[t]\right),
\ \ \mbox{generated by}\ \ 
v_\xi = \bigwedge_{i=1}^r \bigwedge_{j=1}^{\xi_i} v_i \otimes t^{-j}.$$

\begin{theorem} \cite{FL3}
There is a filtration on $\gl_r \otimes \CC\left< X,Y\right>$ such that the adjoint
graded algebra is  isomorphic to $\gl_r \otimes \CC[x^1,x^2]$  and $\grc V_\xi$ is a quotient of
$W^2(\xi)$. Moreover, if $\xi = (n)$ then $\grc V_\xi \cong W^2(\xi)$. \qed
\end{theorem}

\begin{cnj}\label{cnj_2d}
We have $W^2(\xi)^{[k]}\cong \grc \left( V_\xi^{[k]}\right)$.
\end{cnj}

Note that for $\xi = (n)$ we know due to Proposition~\ref{prp_grd} 
that in this conjecture the left hand side is a quotient of the right hand side.

\subsection{Relation to Schubert cells} 

Let us enumerate the basis vectors of $V \otimes \CC[t^{-1},t]$ as follows.
Denote $v_i \otimes t^{-j}$ by $u_{rj-i+1}$.

For a partition $\xi= \xi_1 \ge \dots \ge \xi_r \ge 0$ introduce 
$\eta(\xi)$ as the ordered set of numbers $\eta_1 < \eta_2 < \dots < \eta_n$, $n=|\xi|$,
equal to $lr-s$ with $0 \le s<\xi^t_l$, $l=1,2,\dots $, where $\xi^t$ is the transposed partition.
So 
$$v_\xi = \bigwedge_{i \in \eta(\xi)} u_i.$$

Then our modules are related to Example~\ref{ex_gr} as follows.

\begin{prp} We have
$V^{[k]}_\xi \cong \Gamma(\Sh_{\eta(\xi)}, \oo(k))^*$ and $\prj{V_\xi} \cong Sh_{\eta(\xi)}$.
\end{prp}

\begin{proof}
Note that $V_\xi$ is indeed a submodule of $\bigwedge\nolimits^{|\xi|} U$, where
$$U = V \otimes  \left(t^{-N}\CC[t]/\CC[t]\right), \qquad N = \xi_1,$$
or, in other words, $U = \left<u_i\right>_{i=1\dots Nr}$.

The action of  $\gl_r \otimes \CC\left< X,Y\right>$ on $U$
defines a map from $\gl_r \otimes \CC\left< X,Y\right>$ to 
${\rm End }(U) = \gl_{Nr}$. It is shown in \cite{FL3} that
the image of this map is the ``block upper-triangular'' Lie subalgebra ${\frak p}$
 mapping each $u_{i}$ to a linear combination
of $u_{j}$ with the integer part of $(j-1)/r$ not exceeding the integer part of $(i-1)/r$.
Then $V_\xi^{[k]}$ is the ${\frak p}$--submodule of  
$\bigwedge\nolimits^{|\xi|} U$,
generated by $v_\xi^{\otimes k}$.

Let $\b \subset {\rm End(U)}$ be the Lie algebra of upper-triangular matrices, 
that is endomorphisms mapping each $u_i$ to a linear combination of $u_j$ with $j \le  i$.
Consider the projection  $p^+: {\frak p} \to \b$ along the subspace of strictly lower--triangular matrices.
Since for any $g \in {\frak p}$ 
we have $g v_\xi^{\otimes k} = p^+(g) v_\xi^{\otimes k}$, the subspace  $V_\xi^{[k]}$ is
indeed the $\b$--submodule generated by $v_\xi^{\otimes k}$. 

So we are in the situation of Example~\ref{ex_gr}.
\end{proof}

Set $n=|\xi|$. For a matrix $A = (a_{ij})$, $i=1\dots k$, $j=1\dots n$, of positive integers
 introduce the functional $u^*_A$ on $V^{[k]}_\xi$ by
$$u^*_A = \bigotimes_{i=1}^k \bigwedge_{j=1}^n u_{a_{ij}}^*.$$

\begin{crl}\label{hp}\cite{HP}
Let 
$$\ms_\eta^{[k]} = \{ A = (a_{ij}) \,|\ 
a_{i1} < \dots < a_{in}, \qquad 1 \le a_{1j} \le \dots \le a_{kj} \le \eta_j\}.$$
Then elements $u^*_A$ with $A \in \ms_{\eta(\xi)}^{[k]}$ form a basis in $\left(V^{[k]}_\xi\right)^*$. \qed
\end{crl}

\begin{remark}
Another way to describe the set $\ms_\eta^{[k]}$ is the notion of 
{\em plane partitions}. Recall that a plane partition of shape $\lambda$ is a
filling of the diagram of partition $\lambda$ by non-negative
integers weakly increasing along rows and columns. Let $PP^\lambda(k)$ 
denotes the set of  such plane partitions, filled by integers not exceeding 
$k$, and let $pp^\lambda(k)$ denotes their number.

Then we have the following bijection between $PP^\lambda(k)$ and $\ms_\eta^{[k]}$, where
 $\lambda = (\eta_n - n \ge \dots \ge \eta_1 -1)$. For any plane partition 
from $PP^\lambda(k)$ let $a_{ij}-j$ be the number of
integers in $n-j+1$--th row less than $i$. Then the set $(a_{ij})$ belongs to
$\ms_\eta^{[k]}$ and one can see that it gives us a bijection.
\end{remark}

\begin{remark}
Note that Conjecture~\ref{cnj_2d} implies that $\prj{W^2(\xi)}$ is a degeneration of 
the Schubert variety $\Sh_{\eta(\xi)}$.
\end{remark}

\subsection{Module structure}

\def \fr {{\mathcal F}}

First let us calculate the character of $V_\xi^{[k]}$ as $\gl_r$--module. 
By Corollary~\ref{hp} it can be written as follows.

\begin{prp}
We have 
$$\ch V_\xi^{[k]} = \sum_{A \in \ms_{\eta(\xi)}^{[k]}} \prod_{i,j} x_{a_{ij}\, {\rm mod}\, r},$$
where $a\, {\rm mod}\, r$ takes values from $1$ to $r$. \qed
\end{prp}

Let $\fr_n^r$ be the map from the Grothendiek ring of representations of the symmetric
group $\Sigma_n$ to the Grothendiek ring of $\gl_r$--modules defined by
$$\fr_n^r (\pi) = \left(V^{\otimes n} \otimes \pi\right)^{\Sigma_n}.$$
In other word, it maps the representation of $\Sigma_n$ corresponding to a partition
$\xi$ to the  $\gl_r$--module  corresponding to the same partition if $\xi_{r+1}=0$ and
to zero otherwise.

Recall the notion of {\em skew Schur function}:
$$s_{\lambda\setminus \mu} = 
\sum_{{k_{ij} \ge 1, \, (i,j) \in \lambda \setminus \mu}\atop 
{k_{ij} \le k_{i+1,j}; \, k_{ij} < k_{i,j+1}}}
\prod_{(i,j)\in \lambda \setminus \mu} x_{k_{ij}}.$$
and the representation  of $\Sigma_{|\lambda| - |\mu|}$ 
$$\pi_{\lambda\setminus \mu} = \sum C^\nu_{\lambda \mu} \pi_\nu,$$
where $C^\nu_{\lambda \mu}$ are the Littlewood--Richardson coefficients.
The representation $\pi_{\lambda\setminus \mu}$ corresponds to
the symmetric function $s_{\lambda\setminus \mu}$ as follows.

\begin{prp}\cite{Mac}
 We have
$$\ch \fr_{|\lambda|-|\mu|}^r (\pi_{\lambda\setminus \mu}) =
s_{\lambda\setminus \mu}(x_1,\dots, x_r,0,0,\dots).$$ \qed
\end{prp}

At last for representations $\pi_1$ of $\Sigma_{n_1}$ and $\pi_2$ of $\Sigma_{n_2}$
introduce the outer product 
$$\pi_1 \odot \pi_2 = {\rm Ind}_{\Sigma_{n_1}\times \Sigma_{n_2}}^{\Sigma_{n_1+n_2}}
\pi_1 \boxtimes \pi_2.$$
Note that we have $\fr_n^r (\pi_1 \odot \pi_2) = \fr_n^r (\pi_1) \otimes \fr_n^r (\pi_2)$.

\def \cpf {{\mathbb C}{\rm PF_{Sign}^{[k]}}}

For our purpose introduce the following representation.

\begin{dfn} For a partition $\xi$ set $n = |\xi|$. Then
the {\em Higher Parking Functions} representation  of $\Sigma_{kn}$ is given by
$$\cpf(\xi) = 
\bigoplus_{{\emptyset = \lambda^0 \subset  \lambda^1 \subset \dots \subset \lambda^n = k^n}
\atop{\lambda^s \supset k^{\xi_1^t + \dots +\xi_s^t}, s=1\dots n}} 
\bigodot_{s=1}^n \pi_{\lambda^s\setminus \lambda^{s-1}},
$$ 
where $\xi^t$ is the transposed partition and $k^m = (k \ge \dots \ge k)$, where $m$ is the number
of entrees.
\end{dfn}

\begin{remark}
For $k=1$ this representation is the tensor product of the sign representation
and the representation in $\rho$--parking functions introduced
in \cite{PSt} for $\rho = (1^{\xi_1^t}2^{\xi_2^t} \dots )$, 
that is the partition where each $j$ appears $\xi_j^t$ times, 
and $\xi^t$ is the transposed partition.
\end{remark}

\begin{prp}
We have $V_\xi^{[k]} \cong \fr^r_{kn} \left(\cpf(\xi)\right)$ as $\gl_r$-modules.
\end{prp}

\begin{proof}
It is enough to compare the characters.

For $\Lambda = (\emptyset = \lambda^0 \subset 
\lambda^1 \subset \lambda^2 \subset \dots \subset \lambda^n = k^n)$ 
introduce the set $\ms(\Lambda)$  consisting of matrices $(a_{ij})$ such that 
$$a_{i1} < \dots < a_{in}, \qquad a_{1j} \le \dots \le a_{kj}, \qquad
r(s-1) < a_{ij} \le rs\ \ \mbox{for}\ \ (i,j) \in \lambda^s\setminus \lambda^{s-1}.$$
Then $\ms_{\eta(\xi)}^{[k]}$ is union of $\ms(\Lambda)$ for all 
$\Lambda$ satisfying $\lambda^s \supset k^{\xi_1^t + \dots +\xi_s^t}$, $s=1\dots n$.

Next note that
$$\sum_{A \in \ms(\Lambda)} \prod_{i,j} x_{a_{ij}\, {\rm mod}\, r} = 
\prod_{s=1}^n \sum_{{r\ge k_{ij} \ge 1, \, (i,j) \in \lambda^s \setminus \lambda^{s-1}}\atop 
{k_{ij} \le k_{i+1,j}; \, k_{ij} < k_{i,j+1}}}
\prod_{(i,j)\in \lambda^s \setminus \lambda^{s-1}} x_{k_{ij}}$$
by setting $k_{ij} =  a_{ij}\, {\rm mod}\, r$. And the right hand side is equal to
$$
\prod_{s =1}^n s_{\lambda^s\setminus \lambda^{s-1}}(x_1, \dots, x_r,0,0,\dots) =
\ch \fr^r_{rn} \left(\bigodot_{s=1}^n \pi_{\lambda^s\setminus \lambda^{s-1}}\right)
 $$
\end{proof}

\begin{remark}
The module $W^2(\xi)^{[k]}$ is bi-graded and one of these grading remains in $V_\xi^{[k]}$
as the grading by degree of $t$. This grading can be also viewed in $\cpf(\xi)$
by fixing $|\lambda^1|+ \dots + |\lambda^n|$ in the direct sum.
\end{remark}

\subsection{Dimension formula}

The dimension of $\Gamma(\Sh_{\eta}, \oo(k))$ is given by the 
following {\em Hodge postulation formula}.

\begin{theorem}\cite{HP}\cite{St}
For $\eta =  (\eta_1 < \dots <\eta_n)$ we have
$$ \dim \Gamma(\Sh_{\eta}, \oo(k)) = {\rm Det} \left(
\bin{\eta_j +k-j}{i+k-j}
\right)_{1\le i,j\le n}.
$$ \qed
\end{theorem}

For the case $\xi = (n)$, that is $\eta = (r, 2r,\dots, nr)$,
 let us deduce an explicit formula.
Note that for $k=1$ it gives the higher Catalan
number $C_{n}^{(r)}$  (see e.g. \cite{FL3})
$$\dim V_{(n)} = C_{n}^{(r)}:=\frac{1}{n+1}\bin{r(n+1)}{n}.$$

\begin{theorem} We have
\begin{equation} {\label{ppf}}
\dim V^{[k]}_{(n)} = \prod_{j=1}^n \frac{(jr+k-j)!}{(jr-1)!} \cdot
\frac{(kr+jr-1)!}{(kr+jr-j)!} \cdot \frac{(j-1)!}{(k+j-1)!}.
\end{equation}
\end{theorem}

\begin{proof}
Let $d_{i,j}(k,r) =  \bin{rj +k-j}{i+k-j}$. We have 
$\dim V^{[k]}_{(n)} = {\rm Det} \left(d_{i,j}\right)_{1\le i,j\le n}$. 
First let us make the entrees of this matrix polynomials in $k$ and $r$. 
To do it note that
$$d_{i,j}(k,r) = \frac{(jr+k-j)!}{(k-j+n)! (jr-1)!} d_{i,j}'(k,r)$$
for $d_{i,j}'(k,r) = (k-j+i+1)_{n-i} (rj-i+1)_{i-1}$,
where $(m)_i = m (m+1) \dots (m+i-1)$. So
$$\dim V^{[k]}_{(n)} = \Delta(k,r)\prod_{j=1}^n \frac{(jr+k-j)!}{(k-j+n)!(jr-1)!}, \qquad
\Delta(k,r) = {\rm Det} \left(d_{i,j}'(k,r)\right)_{1\le i,j\le n}.$$

Note that $\Delta(k,r)$ is a polynomial in $k$ and $r$ of degree $n(n-1)/2$ in $k$ and
of the same degree in $r$.

Let us show that $\Delta(k,r)=0$ for $r = b/(a+k)$, $1\le b <a\le n$.
In this case for any $j$ we have 
$$\sum_{i=1}^n d_{i,j}' \left( k \bin{a-b-1}{i-b-1}+ i \bin{a-b}{i-b}\right) = 0,$$
so the rows of this matrix are linearly dependent.

As the polynomials $kr + ar -b$, $1\le b <a\le n$ are irreducible and have no common divisors,
$\Delta(k,r)$ is divisible by $\prod\limits_{1\le b <a\le n} (kr + ar -b)$ and therefore
proportional to this polynomial.

Then note that the maximal degree term  in $\Delta(k,r)$ is equal to
\begin{eqnarray*}
{\rm Det} \left( k^{n-i} (rj)^{i-1}\right)_{1\le i,j\le n} &=& (kr)^{n(n-1)/2} \ 
{\rm Det} \left(j^{i-1} \right)_{1\le i,j\le n}  = \\
&=& (kr)^{n(n-1)/2} \prod_{1\le b <a\le n} (a-b).
\end{eqnarray*}
Summarizing, we obtain
$$\dim V^{[k]}_{(n)} = \prod_{1\le b <a\le n} (a-b)(kr + ar -b)
\prod_{j=1}^n \frac{(jr+k-j)!}{(k-j+n)!(jr-1)!}.$$  
At last writing this in the uniform way using
$$\prod_{1\le b <a\le n} (a-b) = \prod_{j=1}^n (j-1)!, \qquad
\prod_{1\le b <a\le n} (kr + ar -b) = \prod_{j=1}^n \frac{(kr+jr-1)!}{(kr+jr-j)!},$$
and $\prod_{j=1}^n(k-j+n)! = \prod_{j=1}^n (k+j-1)!$,
we obtain the statement of the theorem.
\end{proof}

\begin{crl}
We have
$$\dim W^2(n\omega_1)^{[k]} \ge \prod_{j=1}^n \frac{(jr+k-j)!}{(jr-1)!} \cdot
\frac{(kr+jr-1)!}{(kr+jr-j)!} \cdot \frac{(j-1)!}{(k+j-1)!}.$$\qed
\end{crl}

\begin{remark}
To the best of our knowledge the product formula for the number of plane 
partitions $pp^{\lambda}(k)$ in the case 
$\lambda=(nr+p,(n-1)r+p,\cdots,r+p,p),$ has been obtained for the first time by
R.Proctor (unpublished manuscript dated January 1984, but see \cite{Pr}, 
Corollary~4.1, for the case $r=1$). As far as we aware the first published 
proof of the product formula for $pp^{\lambda}(k)$ in question, is due to C.
Krattenthaler \cite{Kr}. We refer the reader to \cite{St2}, p.550, for an
elegant product formula for the number $pp^{\lambda}(k)$, due to R.Proctor,
 as well as for additional historical comments. We include a (new) proof of
Theorem~3.14 since it directly furnishes the product formula ~\eqref{ppf}
which is more suitable for our purposes.
\end{remark}

\begin{remark}
Using the same method we can prove more generally that for 
$\lambda =(n(r-1)+p,(n-1)(r-1)+p,\cdots,r-1+p)$ we have
$$pp^{\lambda}(k)=
{\rm Det} \left(
\bin{rj +p +k-j}{i+k-j} \right)_{1\le i,j\le n}=$$
$$\prod_{j=1}^n \frac{(jr+p+k-j)!}{(jr+p-1)!} \cdot
\frac{(kr+jr+p-1)!}{(kr+jr+p-j)!} \cdot \frac{(j-1)!}{(k+j-1)!}.$$

The observation that the higher Catalan number $C_{n}^{(r)}$ is equal to the
 number  of plane partitions $pp^{{(n(r-1),(n-1)(r-1),\cdots,r-1)}}(1)$ 
can be generalized to the natural bijection
between the set of {\em trapezoidal paths }  \cite{Lo} of type $(n,p,r)$ 
and the set of plane partitions $PP^{((n-1)(r-1)+p,\cdots,r-1+p)}(1)$.
\end{remark}

\end{document}